\documentclass[12pt,leqno]{amsart}
\usepackage{amsmath,amssymb}
\usepackage[dvips]{graphicx,color}
\usepackage{latexsym,amssymb}
\usepackage{mathrsfs}
\usepackage[abbrev]{amsrefs}

\newtheorem{thm}{Theorem}[section]
\newtheorem{cor}[thm]{Corollaly}
\newtheorem{prop}[thm]{Proposition}
\newtheorem{lem}[thm]{Lemma}

\DeclareMathOperator{\sign}{sign}

 \newenvironment{pf}
    {{\noindent \bf Proof. }}{\hfill $\Box$}

\numberwithin{equation}{section}
\numberwithin{thm}{section}

\begin{document}

\begin{center}\large \bf 
The Leibniz rule for the Dirichlet and the Neumann Laplacian
\end{center}

\footnote[0]
{


E-mail: t-iwabuchi@tohoku.ac.jp

}
\vskip5mm

\begin{center}
Tsukasa Iwabuchi 

\vskip2mm

Mathematical Institute, 
Tohoku University\\
Sendai 980-8578 Japan

\end{center}

\vskip5mm

\begin{center}
\begin{minipage}{135mm}
\footnotesize
{\sc Abstract. }
We study the bilinear estimates in the Sobolev spaces 
 with the Dirichlet and the Neumann boundary condition. 
The optimal regularity is revealed to get such estimates 
in the half space case, 
which is related to not only smoothness of functions and but also 
boundary behavior. 
The crucial point for the proof is how to handle 
boundary values of functions and their derivatives. 

\end{minipage}
\end{center}

\section{Introduction}

We study the bilinear estimates of the form 
\[
\| f g \|_{\dot H^s_p} 
\leq C (\| f \|_{\dot H^s_{p_1}} \| g \|_{L^{p_2}} 
       + \| f \|_{L^{p_3}} \| g \|_{ \dot H^s_{p_4}}),
\]
where $s > 0$ and $p, p_j$ ($j = 1,2,3,4$) satisfy 
$1/p = 1/p_1 + 1/p_2 = 1/p_3 + 1/p_4$. 
The domain is the half space 
$\mathbb R^n_+ 
:= \{ x \in \mathbb R^n \, | \, x_n > 0 \}$, 
and $f,g$ satisfy the boundary condition of either Dirichlet or Neumann type.  
Such inequalities for the Besov spaces are also studied.

The basis of the proof of the bilinear estimates is by 
applying the Leibniz rule and the H\"older inequality. 
This argument works 
in the classical Sobolev spaces $W^{k,p} (\Omega)$ ($k = 1,2,\cdots$), 
where $\Omega$ is an arbitrary domain. 
In the case when $\Omega = \mathbb R^n$, 
such estimates for all regularity $s > 0$ is well-known. 
Classical proof of the bilinear estimates for homogeneous spaces 
can be found in papers by Grafakos and Si~\cite{GS-2012}, Tomita~\cite{Tomi-2010}, 
and it is also proved by the commutator estimates called Kato-Ponce's 
inequality (see a paper by Kato and Ponce~\cite{KaPo-1988}). 
We also refer a book by Runst and Sickel~\cite{RuSi_1996} on the detailed 
analysis of multi-linear estimates, 
and a recent paper by Fujiwara, Georgiev and Ozawa~\cite{FGO-2018} 
who treated higher order fractional Leibniz rule. 
However, when one considers fractional Laplacian on domains, 
there arises difficulty due to how to define 
fractional power and how to handle boundary behavior of functions. 
In general domains, 
we refer to a paper~\cite{IMT-preprint3} 
which studies the bilinear estimates in Besov spaces associated with the 
Dirichlet Laplacian with the regularity $0<s<2$ by means of 
the gradient estimates for the heat equation in $L^p$. 
The exterior domainn case is discussed in a paper~\cite{GeTa-2019}. 
We also refer to several papers by Di Nezza, Palatucci and Valdinoci~\cite{DiPaVa-2012}, 
and Tartar~\cite{Tar-2010} for fractional Sobolev spaces on domains.

In this paper we study in 
function spaces associated with the Dirichlet and the Neumann Laplacian 
in the half space. 
The reason of adapting the half space in this paper  
is just for the sake of simplicity to understand the behavior near the boundary 
clearly, and the obtained result would be able to be applied to other domains. 
We will  understand a reasonable regularity 
for obtaining the bilinear estimates by revealing a roll of derivative 
$\partial _{x_n}$ perpendicular to the boundary. 

Let $A_D$, $A_N$ be the Dirichlet Laplacian $-\Delta |_D$, 
the Neumann Laplacian $-\Delta |_N$, respectively. 
We should note that $A_D, A_N$ can be realized as operators on 
$L^2 (\mathbb R^n_+)$ initially, 
they are regarded as ones of Besov spaces and some spaces of 
distributions by utilizing the uniform boundedness of spectral multipliers 
$\varphi (\theta A_D), \varphi (\theta A_N)$ in $L^1 (\mathbb R^n_+)$ 
with respect to $\theta > 0$. 
Furthermore, the fractional power of $A_D, A_N$ can be defined. 
We refer to related papers \cite{DOS-2002,IMT-RMI,IMT-preprint} 
for boundedness of spectral multipliers, 
\cite{IMT-preprint2} for defining Besov spaces, 
and \cite{Iw-2018} for the fractional Laplacian. 

Let us define spaces of test function spaces, 
Sobolev spaces and Besov spaces following the argument 
\cite{IMT-preprint2} (see also \cite{Tanig} for the Neumann case), 
which are well-defined since 
$e^{-tA_D}$ and $e^{-tA_N}$ satisfies the Gaussian upper bounds. 
The important point there is how to define test function spaces, which 
can give theory of function spaces. 
We take $\phi_0(\cdot) \in C^\infty_0(\mathbb R)$ a non-negative function on $\mathbb R$ such that  
\begin{equation}
\label{325-1}
{\rm supp \, } \phi _0
\subset \{ \, \lambda \in \mathbb R \, | \, 2^{-1} \leq \lambda \leq 2 \, \}, 
\quad \sum _{ j \in \mathbb Z} \phi_0 ( 2^{-j}\lambda) 
 = 1 
 \quad \text{for } \lambda > 0, 
\end{equation}
and $ \{ \phi_j \}_{j \in \mathbb Z}$ is defined by letting 
\begin{gather} \label{325-2}
\phi_j (\lambda) := \phi_0 (2^{-j} \lambda) 
 \quad \text{for }  \lambda \in \mathbb R . 
\end{gather}
Let $\psi$ be a non-negative function such that 
$$
\psi \in C_0 ^\infty (\mathbb R), 
\quad 
\psi (\lambda ) + \sum _{ j \in \mathbb N} \phi _j (\lambda) = 1 
\quad \text{for any } \lambda \geq 0 . 
$$

\vskip3mm 

\noindent 
{\bf Definition (Test function spaces and distributions). } 
{\it 
Let $A = A_D$ or $A_N$. 
\\
{\rm (i)} 
{\rm (}Linear topological spaces
$\mathcal X (A)$ and $\mathcal X' (A)${\rm )} \,
$\mathcal X (A)$ is defined by 
\begin{equation}\notag 
\mathcal X (A) 
:= \big\{ f \in  L^1 (\mathbb R^n_+) \cap \mathcal D (A) 
 \, \big| \, 
    A^{M} f \in L^1(\mathbb R^n_+ ) \cap \mathcal D (A) \text{ for all } M \in \mathbb N 
   \big\} 
\end{equation} 
equipped with the family of semi-norms $\{ p_{A,M} (\cdot) \}_{ M = 1 } ^\infty$ 
given by 
\begin{equation}\notag 
p_{A,M}(f) := 
\| f \|_{ L^1(\mathbb R^n_+)} 
+ \sup _{j \in \mathbb N} 2^{Mj} 
  \| \phi_j (\sqrt{A}) f \|_{ L^1(\mathbb R^n_+)} , 
\end{equation}
and $\mathcal X'(A)$ denotes the topological dual of $\mathcal X(A)$. 

\noindent 
{\rm (ii)}
{\rm (}Linear topological spaces
$\mathcal Z (A)$ and $\mathcal Z' (A)${\rm )} \, 
$\mathcal Z(A)$ is defined by 
\begin{equation}\notag 
\mathcal Z (A)
:= \Big\{ f \in \mathcal X (A) 
 \, \Big| \, 
  \sup_{j \leq 0} 2^{ M |j|} 
    \big\| \phi_j \big(\sqrt{ A } \big ) f \big \|_{L^1(\mathbb R^n_+)} < \infty 
  \text{ for all } M \in \mathbb N
   \Big\} 
\end{equation}
equipped 
with the family of semi-norms $\{ q_{A,M} (\cdot) \}_{ M = 1}^\infty$ 
given by 
\begin{equation}\notag 
q_{A,M}(f) := 
\| f \|_{L^1 (\mathbb R^n_+) }
+ \sup_{j \in \mathbb Z} 2^{M|j|} \| \phi_j (\sqrt{A}) f \|_{L^1(\mathbb R^n_+)},
\end{equation}
and $\mathcal Z'(A)$ denotes the topological dual of $\mathcal Z(A)$. 
}   

\vskip3mm 

\noindent 
{\bf Definition (Besov spaces). } 
{\it 
Let $A = A_D $ or $ A_N$, $s \in \mathbb R$ and $1 \leq p,q \leq \infty$. 
\\
{\rm(i)} $B^s_{p,q} (A) $ is defined by 
\begin{equation}\notag
B^s_{p,q} (A)
:= \{ f \in \mathcal X'(A)
     \, | \, 
     \| f \|_{B^s_{p,q} (A)} < \infty
   \} , 
\end{equation}
where
\begin{equation}\notag 
     \| f \|_{B^s_{p,q} (A)} 
       := \| \psi ( \sqrt{A}) f \|_{L^p  } 
          + \big\| \big\{ 2^{sj} \| \phi_j (\sqrt{A}) f \|_{L^p (\mathbb R^n_+) }  
                 \big\}_{j \in \mathbb  N}
          \big\|_{\ell^q (\mathbb N)}. 
\end{equation}
{\rm(ii)} 
$\dot B^s_{p,q} (A) $ is defined by 
\begin{equation} \notag 
\dot B^s_{p,q}  (A)
:= \{ f \in \mathcal Z'(A)
     \, | \, 
     \| f \|_{\dot B^s_{p,q}(A)} < \infty 
   \} , 
\end{equation}
where
\begin{equation}\notag 
     \| f \|_{\dot B^s_{p,q}(A)} 
       := \big\| \big\{ 2^{sj} \| \phi_j (\sqrt{A}) f\|_{L^p(\mathbb R^n_+)}
                 \big\}_{j \in \mathbb Z}
          \big\|_{\ell^q (\mathbb Z)}. 
\end{equation}
}   

We can also define Sobolev spaces, which were not discussed in \cite{IMT-preprint2} 
(see the well-definedness in  section~\ref{sec_Sobolev}). 

\vskip3mm 

\noindent 
{\bf Definition. } 
{\it 
Let $A= A_D$ or $A_N$, $s \in \mathbb R$ and $1  \leq p \leq \infty$. 
\\
{\rm(i)} $H^s_p (A) $ is defined by 
\begin{equation}\notag
H^s_p (A)
:= \{ f \in \mathcal X'(A)
     \, | \, 
     \| f \|_{H^s_p (A)} := \| (1+A)^{s/2} f \|_{L^p (\mathbb R^n_+)} < \infty
   \} . 
\end{equation}
{\rm(ii)} 
$\dot H^s_p (A) $ is defined by 
\begin{equation} \notag 
\dot H^s_p  (A)
:= \{ f \in \mathcal Z'(A)
     \, | \, 
     \| f \|_{\dot H^s_p(A)}:= \| A^{s/2} f \|_{L^p(\mathbb R^n_+)} < \infty 
   \} . 
\end{equation}
}   

We start by studying derivative operators 
of the normal direction on the boundary $\partial \mathbb R^n_+$ 
(see \cite{Iw-2018-2} for the one dimensional case) 
and derivatives of the other directions. 

\vskip3mm 

\noindent 
{\bf Definition (Derivatives in the sense of distributions). } 
{\it 
\begin{enumerate}
\item[\rm (i)]  
For any 
$f \in \mathcal X'(A_D)$, 
we define $\partial_{x_n} f$ as an element of $\mathcal X'(A_N)$ by 
\begin{equation}\notag 
_{\mathcal X'( A_N)}\langle \partial _{x_n} f , g 
\rangle _{\mathcal X( A_N)} 
:= 
- _{\mathcal X'( A_D)}\langle f, \partial _{x_n} g
\rangle _{\mathcal X(A_D)} 
\quad \text{for any } 
g \in \mathcal X(A_N). 
\end{equation}
For any 
$f \in \mathcal Z'(A_D)$, 
we define $\partial_{x_n} f$ as an element of $\mathcal Z'(A_N)$ by 
\begin{equation}\notag 
_{\mathcal Z'( A_N)}\langle \partial _{x_n} f , g 
\rangle _{\mathcal Z( A_N)} 
:= 
- _{\mathcal Z'( A_D)}\langle f, \partial _{x_n} g
\rangle _{\mathcal Z(A_D)} 
\quad \text{for any } 
g \in \mathcal Z(A_N). 
\end{equation}

\item[\rm (ii)]
For any 
$f \in \mathcal X'( A_N)$, 
we define $\partial_{x_n} f$ as an element of $\mathcal X'(A_D)$ by 
\begin{equation}\notag 
_{\mathcal X'( A_D)}\langle \partial _{x_n} f , g 
\rangle _{\mathcal X( A_D)} 
:= 
- _{\mathcal X'( A_N)}\langle f, \partial _{x_n} g
\rangle _{\mathcal X(A_N)} 
\quad \text{for any } 
g \in \mathcal X( A_D). 
\end{equation}
For any 
$f \in \mathcal Z'( A_N)$, 
we define $\partial_{x_n} f$ as an element of $\mathcal Z'(A_D)$ by 
\begin{equation}\notag 
_{\mathcal Z'( A_D)}\langle \partial _{x_n} f , g 
\rangle _{\mathcal Z( A_D)} 
:= 
- _{\mathcal Z'( A_N)}\langle f, \partial _{x_n} g
\rangle _{\mathcal Z(A_N)} 
\quad \text{for any } 
g \in \mathcal Z( A_D). 
\end{equation}

\item[(iii)] Let $k = 1,2, \cdots, n-1$, 
$A = A_D$ or $A_N$, $X = \mathcal X(A)$ or $\mathcal Z(A)$. 
For $f \in X'$, we define $\partial _{x_k} f$ as an element of $X$ by 
\begin{equation}\notag 
_{X'}\langle \partial _{x_k} f , g 
\rangle _{X} 
:= 
- _{X'} \langle f, \partial _{x_k} g
\rangle _{X} 
\quad \text{for any } 
g \in X. 
\end{equation}

\end{enumerate}
}

\begin{thm}\label{thm_derivative}
Let $s \in \mathbb R$ and $1 \leq p,q \leq \infty$. 
\begin{enumerate}
\item[(i)]
$\partial _{x_n} $ are continuous operators 
from the spaces with the Dirichlet condition 
$\mathcal X(A_D)$, 
$\mathcal X'(A_D)$, 
$\mathcal Z(A_D)$, 
$\mathcal Z'(A_D)$
to those with the Neumann condition 
$\mathcal X(A_N)$, 
$\mathcal X'(A_N)$, 
$\mathcal Z(A_N)$, 
$\mathcal Z'(A_N)$, 
respectively.  

\item[(ii)] 
$\partial _{x_n}$ defines a continuous linear operator 
from $\dot H^s_p (A_D)$ to $\dot H^{s-1}_p (A_N)$ and 
\[
\| \partial _{x_n} f \|_{\dot H^{s-1}_{p}(A_N)} 
\leq C \| f \|_{\dot H^s_{p} (A_D)}, 
\quad 1 < p < \infty.
\]
The same assertion holds for 
the Besov spaces $\dot B^s_{p,q} (A_D)$, $\dot B^s_{p,q}(A_N)$ 
instead of $\dot H^s_p (A_D)$, $\dot H^{s-1}_p (A_N)$, respectively, 
and for the Sobolev and the Besov spaces of inhomogeneous type, 
where $1 < p < \infty$ for the Sobolev spaces. 

\item[(iii)] 
The above assertions {\rm (i)} and {\rm (ii)} also hold by replacing 
$A_D$ and $A_N$ with each other. 

\item[(iv)] 
Let $k = 1,2,\cdots,n-1$, $A = A_D$ or $A_N$. 
Then derivative operators $\partial _{x_k} $ are continuous operators 
from the spaces 
$\mathcal X(A)$, 
$\mathcal X'(A)$, 
$\mathcal Z(A)$, 
$\mathcal Z'(A)$ 
to themselves 
and 
\[
\| \partial _{x_k} f \|_{\dot H^s_p(A)} 
\leq C \| f \|_{\dot H^s_p (A)},
\quad 
\| \partial _{x_k} f \|_{\dot B^{s-1}_{p,q}(A)} 
\leq C \| f \|_{\dot B^s_{p,q} (A)},
\]
where $1 < p < \infty$ for the Sobolev spaces. 
The same assertion holds for the spaces of inhomogeneous type. 
\end{enumerate}

\end{thm}
 
\vskip3mm 

By the above theorem, 
one can understand that $\partial _{x_n}$ changes boundary condition of 
functions essentially 
while the others $\partial _{x_k}$ ($k = 1,2, \cdots, n-1$) do not.

Let us turn to the bilinear estimates. 
Before stating results, we mention a problem to get higher regularity 
of products of functions 
satisfying the Dirichlet and the Neumann boundary condition. 
If the Dirichlet Laplacian acts on a product $fg$ for $f,g$ having 
the Dirichlet boundary condition, one has 
\[
A_D(fg) 
= (A_D f) g - \nabla f \cdot \nabla g + f (A_D g),
\]
and the first and the third term also satisfy the Dirichlet condition 
but $\nabla f, \nabla g$ should have non-zero value on the boundary in general. 
Hence the regularity $s = 2$ case contains an important point, and 
such problem can be found in the Neumann case. 
However, we will have a restriction of regularity only for the Dirichlet case 
and the estimates without restriction for the Neumann case. 
The following is our main theorem.

\begin{thm}\label{thm_bilinear}
Suppose that $p, p_1 , p_2 , p_3, p_4$ satisfy 
\[
1 < p, p_1 ,p_4 < \infty, 
\quad 1 < p_2 , p_3 \leq \infty, 
\quad 
\frac{1}{p} = \frac{1}{p_1} + \frac{1}{p_2} = \frac{1}{p_3} + \frac{1}{p_4}.  
\]
{\rm (i)} 
{\rm (}Dirichlet case{\rm)}
Let $A = A_D, 0 < s < 2 + 1/p$. 
Then there exists $C >0$ such that for any 
$f \in \dot H^s_{p_1} (A_D) \cap L^{p_3} (\mathbb R^n_+)$, 
$g \in L^{p_2} (\mathbb R^n_+) \cap \dot H^s_{p_4}(A_D)$
\begin{equation}\label{1014-1}
\| f g \|_{\dot H^s_p (A_D)} 
\leq C ( \| f \|_{\dot H^s_{p_1}(A_D)} \| g \|_{L^{p_2}} 
  + \| f \|_{L^{p_3}} \| g \|_{\dot H^s_{p_4}(A_D)}
  ).
\end{equation}

\noindent 
{\rm (ii)} 
{\rm (}Neumann case{\rm )}
Let $A = A_N, s > 0$. 
Then there exists $C >0$ such that for any 
$f \in \dot H^s_{p_1} (A_N) \cap L^{p_3} (\mathbb R^n_+)$, 
$g \in L^{p_2} (\mathbb R^n_+) \cap \dot H^s_{p_4}(A_N)$
\begin{equation}\label{1014-2}
\| f g \|_{\dot H^s_p (A_N)} 
\leq C ( \| f \|_{\dot H^s_{p_1}(A_N)} \| g \|_{L^{p_2}} 
  + \| f \|_{L^{p_3}} \| g \|_{\dot H^s_{p_4}(A_N)}
  ).
\end{equation}

\noindent 
{\rm (iii)} 
The corresponding assertion to {\rm (i)} and {\rm (ii)} 
in the inhomogeneous Sobolev spaces hold.
\end{thm}

\begin{thm}\label{thm_counter_example}
Suppose that $s \geq 2 + 1/p$. Then the bilinear estimate \eqref{1014-1} 
of the Dirichlet case does not hold. 
\end{thm}

The result in the Besov spaces also holds. 

\begin{thm}\label{thm_bilinear_2}
Suppose that $p,p_1,p_2,p_3,p_4,q$ satisfy 
\[
1 \leq p, p_1 , p_2 , p_3, p_4  \leq \infty, 
\quad 
\frac{1}{p} = \frac{1}{p_1} + \frac{1}{p_2} = \frac{1}{p_3} + \frac{1}{p_4}.  
\]
Let $s$ be as in Theorem~\ref{thm_bilinear}. 
Then the corresponding bilinear estimates in 
$\dot B^s_{p,q} (A_D)$, $B^s_{p,q}(A_D)$, $\dot B^s_{p,q}(A_N)$, 
$B^s_{p,q} (A_N)$ hold, respectively, 
by replacing the Sobolev spaces with the Besov spaces 
which have the interpolation index~$q$. 
Furthermore, 
if $s > 2+1/p$ or $s = 2+1/p$ with $1 \leq q < \infty$, 
the bilinear estimate does not hold 
for the Dirichlet case. 
\end{thm}

Let us mention multi-linear case. 
There is no restriction of the regularity $s$ for the Neumann case 
which leads to estimates for products of any number of functions. 
On the other hand, $s = 2+1/p$ is optimal for the Dirichlet case. 
Nevertheless, we can show a positive result of some of multi-linear 
estimates for the Dirichlet case. 
Let us state a result for a trilinear inequality as a simplest case.

\begin{cor}\label{cor_tri}
Let $s > 0$, 
$p, p_j $ {\rm(}$j = 1,2, \cdots, 9${\rm)} be such that 
\[
1 < p, p_j < \infty \text{ for } j=1,5,9, 
\quad 1 < p_j \leq \infty \text{ for } j=2,3,4,6,7,8, 
\]
\[
\frac{1}{p} = \frac{1}{p_1} + \frac{1}{p_2} + \frac{1}{p_3}
= \frac{1}{p_4} + \frac{1}{p_5} + \frac{1}{p_6}
= \frac{1}{p_7} + \frac{1}{p_8} + \frac{1}{p_9}. 
\]
Then there exists $C>0$ such that 
\[
\begin{split}
& 
\| fgh \|_{\dot H^s_p (A_D)} 
\\
\leq & 
C ( \| f \|_{\dot H^s_{p_1}(A_D)} \| g \|_{L^{p_2}} \| h \|_{L^{p_3}} 
  + \| f \|_{L^{p_4}} \| g \|_{\dot H^s_{p_5}(A_D)} \| h \|_{L^{p_6}} 
  + \| f \|_{L^{p_7}} \| g \|_{L^{p_8}} \| h \|_{\dot H^s_{p_9}(A_D)} 
  ) . 
\end{split}
\]

\end{cor}

\noindent 
{\bf Remark. } 
One can understand from the proof of Corollary~\ref{cor_tri} (see also \eqref{1101-1}) 
that the multi-linear estimates hold for the product of functions of odd numbers 
but restriction of the regularity appears for the product of even numbers.

\bigskip

Let us give comments about that 
behavior of functions away from the boundary is handled 
similarly to the case $\mathbb R^n$, 
but the main subject is around boundary. 
The cruicial point for the Dirichlet case is: 
The regularity $\alpha = 1/p$ is critical so that 
functions $\nabla f \cdot \nabla g $ for $f,g$ satisfying 
the Dirichlet condition belong to $ H^{\alpha}_p (A_D)$. 
We also notice that $\alpha = 1/p$ is related to considering 
retractions (see page 220 in a book by Triebel~\cite{Triebel_1978}). 
This applied to $A_D (fg) $ leads to reach at the regularity number 
$s = 2 + 1/p$ in Theorem~\ref{thm_bilinear}. 
It is characteristic of the two theorems that 
$\nabla f \cdot \nabla g \not \in \dot H^{\frac{1}{p}}_p (A_D)$ 
breaks down the bilinear estimates in Theorem~\ref{thm_bilinear} (i)
for $s = 2+ 1/p$, 
while $C_0^\infty (\mathbb R^n_+)$ is dense in the Sobolev space 
$H^s_p (\mathbb R^n_+)$ with $s \leq 1/p$ 
defined by the restriction of functions on $\mathbb R^n$ to $\mathbb R^n_+$. 
Here we mention a paper by Killip, Visan and Zhang~\cite{KVZ-2016}, 
where the case when $s < 1 + 1/p$ is studied for exterior domains. 
They obtained the bilinear estimates for $s < 1 + 1/p$, showing that 
the equivalence of 
$(-\Delta|_D) f \in L^p $ and $(-\Delta_{\mathbb R^n}) f \in L^p$ for 
$f\in C_0 ^\infty (\Omega)$, where $-\Delta|_D$ is the Dirichlet 
Laplacian on $\Omega$, $-\Delta_{\mathbb R^n}$ is the Laplacian on $\mathbb R^n$. 
Here it would be reasonable to conjecture that: 
{\it $s = 2 + 1/p$ is the universal upper bound 
for the bilinear estimate \eqref{1014-1} for the Dirichlet case 
in any domain. 
}

It would be plausible that the optimality of $s = 2+ 1/p$ is due to 
the high spectral component affecting the local behavior of functions 
around the boundary. 
As for the low spectrum, which is essetial for the homogeneous spaces, 
it depends on domains. 
The bounded domain case has no restriction, but the possible regularity 
in the exterior domain case is restricted to smaller range 
because of the slower decay of gradient estimates for the heat kernel 
(see papers~\cite{GeTa-2019,IshiKabe-2007}). 

\bigskip

In contrast, the situation is quite different for the Neumann 
condition in spite of that each of $\nabla f, \nabla g$ for $f,g$ with the Neumann 
condition can not expected to satisfy again the Neumann condition. 
The reason is due to that $\nabla f$, $\nabla g$ satisfy the 
Dirichlet condition, which give the Neumann condition for the product 
$\nabla f \cdot \nabla g$, 
and hence, we could expect no restriction of the regularity $s$ for the bilinear estimates.

\vskip3mm 

This paper is organized as follows. 
In section 2, we prepare some important estimates and relations between 
two cases of $\mathbb R^n_+$ and $\mathbb R^n$ in the Sobolev and the Besov 
spaces. In section 3, Theorem~\ref{thm_derivative} is proved. 
Section 4 is devoted to proving bilinear and trilinear estimates of 
Theorems~\ref{thm_bilinear}, \ref{thm_bilinear_2} and Corollary~\ref{cor_tri}. 
In section 5, counterexamples in Theorem~\ref{thm_counter_example} 
will be given.


\vskip3mm 

\noindent 
{\bf Notations. }
Upper and lower half spaces are written as $\mathbb R^n_+ 
:= \{ x \in \mathbb R^n \, | \, x_n > 0 \}$, 
$\mathbb R^n_- 
:= \{ x \in \mathbb R^n \, | \, x_n < 0 \}$. 
We often write $x \in \mathbb R^n$ as 
$x = (x' , x_n)$, where 
$ x' \in \mathbb R^{n-1}, x_n \in \mathbb R$. 
The fractional Laplacian in $\mathbb R^n$ is written as 
\[
\Lambda := \mathcal F^{-1} |\xi| \mathcal F .
\]
$\nu$ denotes the outer unit normal vector on the boundary $\partial \mathbb R^n_+$. 
We often omit the domain $\mathbb R^n_+$ in the norm of $L^p (\mathbb R^n_+)$, 
and write $\mathbb R^n$ clearly, more concretely, 
\[
\| f \|_{L^p (\mathbb R^n_+)} = \| f \|_{L^p}, 
\quad 
\| f \|_{L^p (\mathbb R^n)}, 
\quad 
\| f \|_{\dot B^s_{p,q} (\mathbb R^n)}, 
\]
For any function $f$ on $\mathbb R^n_+$, 
let $f_{odd}, f_{even}$ be odd, even extention of $f$ with respect to 
$x_n $ component, respectively, namely,  
\begin{equation}\label{even_odd}
f_{odd} := 
\begin{cases} 
f(x) , & x_n > 0 ,
\\
-f(-x) , & x_n < 0 , 
\end{cases}
\qquad 
f_{even} := 
\begin{cases} 
f(x) , & x_n > 0 ,
\\
f(-x) , & x_n < 0 . 
\end{cases}
\end{equation}

\section{Preliminary}

We prepare useful lemmas to prove our theorems in this section. 
Let us start by enumerating known facts; 
The boundedness of the Riesz transformation in $\mathbb R^n$
(see e.g. a book by Stein~\cite{Stein_1970}),  
the real interpolation of the Sobolev spaces and the Besov spaces
(see~\cite{BL_1976,Iw-2018,Triebel_1983}). 
Then we will state lemmas which are fundamental for our proof. 

\begin{lem}
{\rm (i)} {\rm (}Boundedness of Riesz transform{\rm)} 
Let $1 < p < \infty $. Then a constant $C>0$ exists  such that 
 \[
\| (-\Delta)^{-1/2} \partial _k f \|_{L^p (\mathbb R^n)} 
\leq C \| f \|_{L^p (\mathbb R^n)},
\quad  k = 1,2,\cdots ,n 
\]

\noindent 
{\rm (ii)} {\rm (}Real interpolation{\rm)}
Let $0 < \theta < 1$, $s,s_0,s_1 \in \mathbb R$ and $1 \leq p,q,q_0,q_1 \leq \infty$. 
Assume that $s_0 \not = s_1$ and 
$s = (1-\theta)s_0 + \theta s_1$. Then 
\begin{gather}
\notag 
\big(  \dot B^{s_0}_{p,q_0} , \, \dot B^{s_1}_{p,q_1} \big)_{\theta , q} 
= \dot B^s_{p,q}, 
\quad 
\big(  B^{s_0}_{p,q_0} , \,  B^{s_1}_{p,q_1} \big)_{\theta , q} 
= B^s_{p,q}, 
\\ \notag 
\big(  \dot H^{s_0}_{p}  , \, \dot H^{s_1}_{p} \big)_{\theta , q} 
= \dot B^s_{p,q}, 
\quad 
\big(  H^{s_0}_{p} , \,  H^{s_1}_{p} \big)_{\theta , q} 
= B^s_{p,q}, 
\end{gather}
where 
$\dot B^s_{p,q} := \dot B^s_{p,q} (\mathbb R^n), \dot B^s_{p,q} (A)$, 
$B^s_{p,q} := B^s_{p,q} (\mathbb R^n), B^s_{p,q} (A)$, 
$\dot H^s_{p} := \dot H^s_{p} (\mathbb R^n), \dot H^s_{p} (A)$, 
$ H^s_{p} := H^s_{p} (\mathbb R^n), H^s_{p} (A)$, 
respectively. 

\end{lem}

\begin{lem}\label{lem_1}
Let $s \geq 0, 1 \leq p \leq \infty$ and $f \in L^p(\mathbb R^n_+)$. 
Then,
$A_D ^{s/2} f \in L^p (\mathbb R^n_+)$ if and only if 
$\Lambda ^s f_{odd} \in L^p (\mathbb R^n)$. 
Also, $A_N ^{s/2} f \in L^p (\mathbb R^n_+)$ if and only if 
$\Lambda ^s f_{even} \in L^p (\mathbb R^n)$. Furthermore, 
\[
2^{\frac{1}{p}} \| A_D ^{s/2} f \|_{L^p (\mathbb R^n_+)} 
=  \| \Lambda ^s f_{odd} \|_{L^p (\mathbb R^n)}, 
\quad 
2^{\frac{1}{p}} \| A_N ^{s/2} f \|_{L^p (\mathbb R^n_+)} 
= \| \Lambda ^s f_{even} \|_{L^p (\mathbb R^n)} . 
\]

\end{lem}


\begin{pf}
We consider the Dirichlet case only, since the Neumann case follows analogously 
by using even extention instead of odd one.

We start by proving in the case when $0 < s \leq 2$. 
Let $P_s (t,x) := \mathcal F^{-1} [e^{-t|\xi|^s}](x)$ for 
$t > 0$ and $x \in \mathbb R^n$. 
Since the kernel of $e^{-t A_D ^{s /2}}$ is given by the difference of 
$P_s$, we write 
\begin{equation} \notag 
\begin{split}
e^{-t A_D ^{s /2}}f 
=
& 
\int_{\mathbb R^n_+} \big(P_s (t,x-y) - P_s (t, x'-y' , x_n + y_n)\big) f(y) dy
\\
= 
& \int_{\mathbb R^n} P_s (t,x-y)  f_{odd}(y) dy  
\quad \text{ for } x \in \mathbb R^n _+ , 
\end{split}
\end{equation}
and 
\begin{equation}\label{0512}
\lim _{t \to 0}\frac{e^{-t A_D ^{s /2}}f - f}{t}
= ( \Lambda^{s } f_{odd} ) |_{\mathbb R^n_+} ,
\end{equation}
which implies that 
$A_D ^{s /2 } f \in L^p (\mathbb R^n_+)$ 
and $\Lambda ^s f_{odd} \in L^p (\mathbb R^n)$ are 
equivalent, 
$2^{\frac{1}{p}} \| A_D ^{s/2} f\|_{L^p (\mathbb R^n_+)} = 
 \| \Lambda ^s f_{odd} \|_{L^p (\mathbb R^n)}$.

Let us consider the case when $2 < s \leq 4$. 
For $f \in L^p (\mathbb R^n_+)$ with $A^{s/2} _D f \in L^p (\mathbb R^n_+) $, 
we can see that $A f = -\Delta f \in L^p (\mathbb R^n_+)$, 
$\Delta f_{odd}$ is given by the odd extention of $\Delta f$ and 
$\Delta f_{odd} \in L^p (\mathbb R^n)$, since 
for any $\varphi \in \mathbb R^n$ 
$$
\int_{\mathbb R^n} f_{odd} \Delta \varphi dx 
= \int _{\partial \mathbb R^n_+ \cup \partial \mathbb R^n_-} f_{odd} \nabla \varphi \cdot \nu  dS 
  - \int _{\partial \mathbb R^n_+ \cup \partial \mathbb R^n_-}
   (\nabla f_{odd} \cdot \nu )\varphi  dS
  + \int _{\mathbb R^n} (\Delta f_{odd} ) \varphi dx, 
$$
the first two terms in the right hand side are zero thanks to 
$f_{odd}$ vanishing at $x_n = 0$ and
even property of $\partial_{x_n} f_{odd}$ 
and the integrals of $f_{odd}$ on $\partial \mathbb R^n_+ , \partial \mathbb R^n_-$ 
are justified by the well-definedness of the trace operators 
of $f , \nabla f $ on $ \partial \mathbb R^n_+$ with value in 
$L^p(\partial \mathbb R^n_+)$
for $f $ with $f, \Delta f \in L^p (\mathbb R^n _+)$. 
Hence, $A_D ^{s/2} f \in L^p (\mathbb R^n_+)$ implies that 
\[
\begin{split}
A_D ^{s/2} f 
=
&  A_{D} ^{s/2-2} A_D  f 
= A_D ^{s/2 -2} \big( (-\Delta f_{odd}) |_{\mathbb R^n _+} \big)
= \Big( \Lambda ^{s-2} (-\Delta f_{odd}) \Big) \Big|_{\mathbb R^n _+} 
\\
=
&  ( \Lambda ^{s} f_{odd} ) |_{\mathbb R^n _+} 
\in L^p (\mathbb R^n _+),
\end{split}
\]
which proves $\Lambda ^s f_{odd} \in L^p (\mathbb R^n)$, 
since $\Lambda ^s f_{odd}$ is an odd function.  
Conversely, let $\Lambda ^s f_{odd} \in L^p (\mathbb R^n)$. 
Here $f_{odd}, -\Delta f _{odd} \in L^p (\mathbb R^n)$ implies 
the well-definedness of the trace operator of 
$f_{odd} |_{\mathbb R^n_+}$, 
which implies $\Lambda^2 f_{odd}(x) = A_D f (x)$ 
for almost every $x \in \mathbb R^n_+$ by using the equality 
\eqref{0512}. 
Now, by applying the result in the case when $0 < s \leq 2$ 
proved above to a function $A_D f$, 
we get the equivalence of $A_D ^{s/2-2} (A_D f) \in L^p (\mathbb R^n_+)$ 
and $\Lambda ^{s-2} (-\Delta f_{odd}) \in L^p (\mathbb R^n)$. 
Therefore we have that $\Lambda ^s f_{odd} \in L^p (\mathbb R^n)$ 
gives $A_D ^{s/2} f \in L^p (\mathbb R^n_+)$. 

By the above argument together with the induction, 
we get the result for $k < s \leq k + 2$ for any even number $k$, 
which completes the proof. 
\end{pf}

\begin{lem}
Suppose $1 < p < \infty$, $0 < s < 1/p$. 
Let $\chi _{x_n > 0}$ denote the characteristic function on 
$\{ x\in\mathbb R^n \, | \, x_n > 0 \}$. 
Then there exists 
$C > 0$ such that for any $f \in \dot H^s_p (\mathbb R^n)$ 
\begin{equation}\label{1011-2}
\| \chi _{x_n > 0} f \|_{\dot H^s_p (\mathbb R^n)} 
\leq C \| f \|_{\dot H^s_p (\mathbb R^n)} , 
\end{equation}
\begin{equation}\label{1011-2-0}
\| ( \sign x_n) f \|_{\dot H^s_p (\mathbb R^n)} 
\leq C \| f \|_{\dot H^s_p (\mathbb R^n)} .
\end{equation}
Let $f $ be a function on $\mathbb R^n _+$. Then $f_{odd}$ and $f_{even}$ 
enjoy 
\begin{equation}\label{1011-2-2}
\| f_{odd} \|_{\dot H^s_p (\mathbb R^n)} 
\leq C \| f_{even} \|_{\dot H^s_p (\mathbb R^n)}, 
\quad 
\| f_{even} \|_{\dot H^s_p (\mathbb R^n)} 
\leq C \| f_{odd} \|_{\dot H^s_p (\mathbb R^n)} .
\end{equation}
\end{lem}

\begin{pf}
For $\varphi \in C^\infty (\mathbb R)$ with 
$0 \leq \varphi \leq 1$, 
$\varphi (x_n) = 1$ for $x_n \geq 1$ and 
$\varphi (x_n) = 0$ for $x_n \leq 1/2$, 
put 
\[
\varphi _{\varepsilon} = \varphi _{\varepsilon} (x_n) := 
\varphi (\varepsilon ^{-1} x_n) 
\text{ for any } x_n \in \mathbb R. 
\]
Let us start by proving the uniform boundedness with respect to 
$\varepsilon > 0$, 
\begin{equation}\label{1011-1}
\| \varphi _{\varepsilon} f \|_{\dot H^s_p(\mathbb R^n)} 
\leq C \| f \|_{\dot H^s_p(\mathbb R^n)}. 
\end{equation}
By Bony's paraproduct formula (see \cite{Bo-1981}), 
we consider the frequency decomposition 
\[
\varphi _{\varepsilon} f 
=
\Big( \sum _{k \leq l+3} + \sum _{k > l+3} \Big)
\Big( \phi_k (\sqrt{-\Delta}) \varphi_\varepsilon\Big)
\Big( \phi_l (\sqrt{-\Delta}) f \Big)
=: (\varphi _{\varepsilon}f)_{I} + (\varphi _{\varepsilon} f)_{II}, 
\]
where the first one has component such that 
frequency of $f$ higher than or comparable with that of $\varphi_\varepsilon$, 
and the second one has the other such that 
frequency of $f$ lower than that of $\varphi_\varepsilon$. 
Then applying the bilinear estimate in the Sobolev spaces in 
$\mathbb R^n$ to the first term gives that 
\[
\| (\varphi _{\varepsilon} f)_{I} \|_{\dot H^s_p(\mathbb R^n)} 
\leq C \| \varphi_{\varepsilon} \|_{L^\infty(\mathbb R^n)} 
    \| f \|_{\dot H^s_p(\mathbb R^n)} 
\leq C \| f \|_{\dot H^s_p(\mathbb R^n)} ,
\]
since $f$ has higher frequency than that of $\varphi_\varepsilon $. 
As for the second term, applying the bilinear estimate in the Sobolev spaces 
for the component $x_n$ with indices 
$p_1$ and $p_2$ such that $1/p = 1/p_1 + 1/p_2$, $s = 1/p_1$, $s = 1/p - 1/p_2$  
and the embedding $\dot H^s_p (\mathbb R) \subset L^{p_2} (\mathbb R)$ 
give that 
\[
\begin{split}
\| (\varphi _{\varepsilon} f)_{II} \|_{\dot H^s_p(\mathbb R^n)} 
\leq 
& C \big\| 
       \| \varphi_\varepsilon  \|_{\dot H^s_{p_1}(\mathbb R_{x_n})} 
       \| f \|_{L^{p_2} (\mathbb R_{x_n})}
       \big\| _{L^p (\mathbb R^{n-1}_{x'})}
\\
\leq 
& C        \| \varphi_\varepsilon  \|_{\dot H^{\frac{1}{p_1}}_{p_1}(\mathbb R_{x_n})} 
       \big\| 
       \| f \|_{\dot H^s_{p} (\mathbb R_{x_n})}
       \big\| _{L^p (\mathbb R^{n-1}_{x'})}
\\
\leq 
& C        \| \varphi_1  \|_{\dot H^{\frac{1}{p_1}}_{p_1}(\mathbb R_{x_n})} 
       \big\| \mathcal F^{-1}|\xi_n|^s \mathcal F f \big\|_{L^{p} (\mathbb R^n)}. 
\end{split}
\]
Here it should be noted that when we apply the bilinear estimate above, 
the frequency of 
$(\varphi_\varepsilon f)_{II}$ is restricted to $\xi_n$ direction, 
since 
$\varphi _\varepsilon$ have only the frequency component for $x_n$ 
and its frequency higher than $f$, and 
$s < 1/p$ implies $p_2 < \infty$. By applying 
the Fourier multiplier theorem to a Fourier multiplier 
$|\xi_n|^s / |\xi|^s$, we have 
\begin{equation}\label{1115-1}
\big\| \mathcal F^{-1}|\xi_n|^s \mathcal F f \big\|_{L^p (\mathbb R^n)}
\leq C \| \mathcal F^{-1}|\xi|^s \mathcal F f \|_{L^p(\mathbb R^n)} 
= C \| f \|_{\dot H^s_p(\mathbb R^n)},
\end{equation}
which completes the proof of \eqref{1011-1}. 
Since $\dot H^s_p (\mathbb R^n)$ is a reflexive Banach space 
and $\varphi _ \varepsilon f$ converges to $\chi_{x_n > 0} f$ 
weakly in $\dot H^s_p (\mathbb R^n)$ as $\varepsilon \to 0$, we obtain 
\[
\| \chi_{x_n > 0} f \|_{\dot H^s_p(\mathbb R^n)} 
\leq \liminf _{\varepsilon \to 0} 
\| \varphi_\varepsilon f \|_{\dot H^s_p(\mathbb R^n)}
\leq C \| f \|_{\dot H^s_p(\mathbb R^n)} ,
\]
by taking a subsequence of $\{\varepsilon > 0 \}$ 
if necessary, which proves \eqref{1011-2}. 
The inequality \eqref{1011-2-0} follows from $\sign x_n = 2\chi_{x_n>0} - 1$ 
and \eqref{1011-2}. 
The last inequalities \eqref{1011-2-2} are obtained by 
$f_{even} - f_{odd} = 2 \chi _{x_n >0}$ and 
\eqref{1011-2}. 
\end{pf}

\begin{lem}\label{lem_1102}
Let $1 < p < \infty$, $s \geq 0$, $f \in L^p (\mathbb R^n_+)$. Then 
\begin{equation}\label{1011-3}
\| A_D^{s/2} \partial _{x_n}f \|_{L^p} 
\leq C \| f \|_{\dot H^{s+1}_{p} (A_N)}, 
\quad 
\| A_N^{s/2} \partial _{x_n}f \|_{L^p} 
\leq C \| f \|_{\dot H^{s+1}_{p} (A_D)} 
\end{equation}
provided that the left hand sides are finite, respectively.
Let $1 < p < \infty$, $0 \leq s < 1/p$. Then 
\begin{equation}\label{1011-4}
\| A_D^{s/2} \partial _{x_n}f \|_{L^p} 
\leq C \| f \|_{\dot H^{s+1}_{p} (A_D)}, 
\quad 
\| A_N^{s/2} \partial _{x_n}f \|_{L^p} 
\leq C \| f \|_{\dot H^{s+1}_{p} (A_N)}, 
\end{equation}
provided that the left hand sides are finite, respectively.
\end{lem}

\begin{pf}
We start by proving the first inequality of \eqref{1011-3}. 
Let $f \in L^p (\mathbb R^n_+) \cap \dot H^{s+1}_p(A_N)$, 
which also satisfies 
$f_{even} \in L^p (\mathbb R^n) \cap \dot H^{s+1}_p (\mathbb R^n)$ 
by Lemma~\ref{lem_1}. Firstly, since 
$\Lambda f_{even} \in L^p(\mathbb R^n)$ and the boundedness of 
the Riesz transform give $\nabla f \in L^p (\mathbb R^n_+)$, 
we can see that 
the trace of $f $ in $ L^p (\partial \mathbb R^n_+)$ makes sense 
by the trace theorem (see e.g. \cite{Triebel_1978}). 
Observe 
$(\partial _{x_n} f)_{odd} = \partial _{x_n} f_{even}$ 
which is assured by 
\[
\begin{split}
& 
\int_{\mathbb R^n} (\partial_{x_n} f) _{odd} (x) \varphi (x) dx
\\
= 
& 
\int_{x_n > 0} \partial_{x_n} f (x) \varphi (x) dx
+ \int_{x_n < 0} \partial _{x_n} ( f(x', -x_n)) \varphi(x) dx
\\
=
&
\lim _{\varepsilon \to 0} 
\Big( 
-\int_{\mathbb R^{n-1}} f(x',\varepsilon)dx' 
+\int_{\mathbb R^{n-1}} f(x',\varepsilon)dx'  
\Big)
 - \int_{\mathbb R^n} f_{even} \partial_{x_n}\varphi (x) dx 
\\
=
&
 - \int_{\mathbb R^n} f_{even} \partial_{x_n}\varphi (x) dx , 
\quad \varphi \in \mathcal S(\mathbb R^n). 
\end{split}
\]
Here we should note that the above integrals on $\mathbb R^{n-1}$ is zero, 
since this is justified by the well-definedness 
of the trace operator of $f$ with value in
 $L^p(\partial \mathbb R^n_+)$. 
Lemma~\ref{lem_1} and the boundedness of Riesz transform imply 
\[
\| A_D^{s/2} \partial _{x_n}f \|_{L^p} 
\leq C \| \Lambda^s (\partial_{x_n} f)_{odd} \|_{L^p} 
\leq C \| \Lambda^s \partial_{x_n} (f_{even}) \|_{L^p} 
\leq C \| \Lambda^{s+1} f_{even} \|_{L^p} 
\leq C \| f \|_{\dot H^s_{p} (A_N)} , 
\]
which proves the first inequality of \eqref{1011-3}. 
The second one follows analogously. In fact, 
let $f \in L^p (\mathbb R^n_+) \cap \dot H^{s+1}_p (A_D)$ which also satisfies 
$f_{odd} \in L^p (\mathbb R^n ) \cap \dot H^{s+1}_p (\mathbb R^n)$,  
and the trace of $f$ with value in $L^p (\partial \mathbb R^n_+)$ makes sense. 
Furthermore, the trace of $f$ is zero, since odd function $f_{odd}$ 
is zero on $\{x_n = 0\}$. 
Observe $(\partial_{x_n} f) _{even} = \partial_{x_n} f_{odd}$ which 
is assured by 
\[
\begin{split}
& 
\int_{\mathbb R^n} (\partial_{x_n} f) _{even} (x) \varphi (x) dx
\\
= 
& 
\int_{x_n >0} \partial _{x_n} f(x) \varphi (x) dx 
- \int_{x_n <0} \partial_{x_n} (f(x',-x_n)) dx 
\\
= 
& 
\lim_{\varepsilon \to 0} 
\Big( 
-\int_{\mathbb R^{n-1}} f(x',\varepsilon)dx' 
-\int_{\mathbb R^{n-1}} f(x',\varepsilon)dx'  
\Big) 
- \int_{\mathbb R^n} f_{odd} \partial_{x_n}\varphi (x) dx 
\\
= 
& - \int_{\mathbb R^n} f_{odd} \partial_{x_n}\varphi (x) dx , 
\quad \varphi \in \mathcal S(\mathbb R^n), 
\end{split}
\]
where the integrals on $\mathbb R^{n-1}$ vanishes thanks to 
the trace of $ f$ is zero. Therefore, we obtain 
\[
\| A_N^{s/2} \partial _{x_n}f \|_{L^p} 
\leq C \| \Lambda^s (\partial_{x_n} f)_{even} \|_{L^p} 
\leq C \| \Lambda^s \partial_{x_n} (f_{odd}) \|_{L^p} 
\leq C \| \Lambda^{s+1} f_{odd} \|_{L^p} 
\leq C \| f \|_{\dot H^s_{p} (A_D)} ,
\]
which proves the first inequality of \eqref{1011-3}.

We turn to prove the second one \eqref{1011-4}. 
It follows from \eqref{1011-2-2} that 
\[
\| \Lambda^s F_{odd} \|_{L^p(\mathbb R^n)} 
\leq C \| \Lambda ^s F_{even} \|_{L^p(\mathbb R^n)}, 
\quad 
\| \Lambda^s F_{even} \|_{L^p(\mathbb R^n)} 
\leq C \| \Lambda ^s F_{odd} \|_{L^p(\mathbb R^n)}. 
\]
These inequalities for $F = \partial _{x_n} f$ and the similar 
argument to prove \eqref{1011-3} give that 
\[
\| A_D^{s/2} \partial _{x_n}f \|_{L^p} 
\leq C \| \Lambda^s (\partial_{x_n} f)_{even} \|_{L^p(\mathbb R^n)} 
\leq C \| (\Lambda^s \partial_{x_n}) f_{odd} \|_{L^p(\mathbb R^n)} 
\leq C \| f \|_{\dot H^s_p (A_D)} , 
\]
\[
\| A_N^{s/2} \partial _{x_n}f \|_{L^p} 
\leq C \| \Lambda^s (\partial_{x_n} f)_{odd} \|_{L^p(\mathbb R^n)} 
\leq C \| (\Lambda^s \partial_{x_n}) f_{even} \|_{L^p(\mathbb R^n)} 
\leq C \| f \|_{\dot H^s_p (A_N)} , 
\]
which proves \eqref{1011-4}. 
\end{pf}

\begin{lem}\label{lem_1102-2} 
Let $s \geq 0, 1 < p < \infty$ and $f \in L^p(\mathbb R^n_+)$. 
Then
\[
\| A_D ^{s/2} \partial _{x_k} f \|_{L^p (\mathbb R^n_+)} 
\leq C \| \Lambda ^{s+1} f_{odd} \|_{L^p (\mathbb R^n)}, 
\quad 
\| A_N ^{s/2} \partial _{x_k} f \|_{L^p (\mathbb R^n_+)} 
\leq C  \| \Lambda ^{s+1} f_{even} \|_{L^p (\mathbb R^n)}  
\]
for $k = 1,2, \cdots, n-1$.

\end{lem}

\begin{pf}
Let us prove the first inequality. 
Let $f \in L^p (\mathbb R^n_+)$ be such that 
$\Lambda^{s+1} f_{odd} \in L^p (\mathbb R^n_+)$. By Lemma~\ref{lem_1} and 
the boundedness of the Riesz transform, 
\[
\| A_D^{s/2} \partial _{x_k} f \|_{L^p} 
= \frac{1}{2^{\frac{1}{p}}}
  \| \Lambda ^s (\partial_{x_k} f)_{odd} \|_{L^p (\mathbb R^n)}
= \frac{1}{2^{\frac{1}{p}}} 
  \| \Lambda ^{s} \partial _{x_k} f_{odd} \|_{L^p (\mathbb R^n)} 
\leq C \| \Lambda ^{s+1} f_{odd} \|_{L^p (\mathbb R^n)}.
\]
The second inequality follows analogously. 
\end{pf}

\begin{lem}\label{lem_1115-1}
Let $s \in \mathbb R$, $1 \leq p,q \leq \infty$. Then 
\[
\| f \|_{\dot B^s_{p,q} (A_D)} \simeq \| f_{odd} \|_{\dot B^s_{p,q} (\mathbb R^n)}, 
\quad 
\| f \|_{\dot B^s_{p,q} (A_N)} \simeq \| f_{even} \|_{\dot B^s_{p,q} (\mathbb R^n)},
\]
\[
\| \partial_{x_n}f \|_{\dot B^s_{p,q} (A_D)} 
\leq C \| f \|_{\dot B^{s+1}_{p,q} (A_N)} , 
\quad 
\| \partial_{x_n}f \|_{\dot B^s_{p,q} (A_N)} 
\leq C \| f \|_{\dot B^{s+1}_{p,q} (A_D)} ,
\]
\[
\| \partial_{x_k}f \|_{\dot B^s_{p,q} (A)} 
\leq C \| f \|_{\dot B^{s+1}_{p,q} (A)} 
\quad \text{ for } A = A_D, A_N, \,\,\,\,  k = 1,2,\cdots, n-1.
\]
The corresponding equivalence and inequalities for the inhomogeneous spaces 
$B^s_{p,q}$ also hold. 
\end{lem}

\begin{pf}
Let $M \in \mathbb N$ be such that $M > s/2$. It follows from Theorem~1.3 
in \cite{Iw-2018} that 
\[
\| f \|_{\dot B^s_{p,q} (A_D)} 
\simeq 
  \left\{ \int_0^\infty 
  \Big( t^{-\frac{s}{2}} \| (tA_D)^M e^{-tA_D}f \|_{L^p} \Big)^q 
  \right\}^{\frac{1}{q}}.
\]
Observing that 
\[
(tA_D)^M e^{-tA_D}f (x)
= 
\int_{\mathbb R^n} 
 (-t\Delta ^M G_t) (x-y) f_{odd}(y) dy, 
\]
where $G_t (x) := (4\pi t)^{-\frac{n}{2} } e^{-\frac{|x|^2}{4t}}$, 
we get 
\[
\| f \|_{\dot B^s_{p,q} (A_D)} 
\simeq 
  \left\{ \int_0^\infty 
  \Big( t^{-\frac{s}{2}} \| (-t\Delta)^M e^{t \Delta}f_{odd} \|_{L^p (\mathbb R^n)} \Big)^q 
  \right\}^{\frac{1}{q}}
\simeq 
\| f_{odd} \|_{\dot B^s_{p,q} (\mathbb R^n)},  
\]
which proves the Dirichlet Laplacian case of the homogeneous type. 
The Neumann case follows analogously by means of even extension instead of 
odd one. 
The inhomogeous case is proved by a similar argument to the above and 
using equivalent norms of Besov spaces by semigroup  
(see Theorem~7.2 in \cite{Iw-2018})
\[
\| f \|_{B^s_{p,q}(A)} 
\simeq \| \psi (A) f\|_{L^p}
+  \left\{ \int_0^1
  \Big( t^{-\frac{s}{2}} \| (tA_D)^M e^{-tA_D}f \|_{L^p} \Big)^q 
  \right\}^{\frac{1}{q}}.
\]
We have obtained the norm equivalence. 

We turn to prove the inequalities for $\partial _{x_n} f$. 
Following the proof of \eqref{1011-3} and applying the equivalence obtained above, 
we see that 
\[
\begin{split}
\| \partial _{x_n} f \|_{\dot B^s_{p,q} (A_D)}
\leq 
& 
C 
\| ( \partial _{x_n} f )_{odd}\|_{\dot B^s_{p,q} (\mathbb R^n)}
\leq C 
\| \partial _{x_n} f_{even}\|_{\dot B^s_{p,q} (\mathbb R^n)}
\leq C 
\| f_{even}\|_{\dot B^{s+1}_{p,q} (\mathbb R^n)}
\\
\leq 
& 
C 
\| f\|_{\dot B^{s+1}_{p,q} (A_N)},
\end{split}
\]
and similarly, 
\[
\| \partial _{x_n} f \|_{\dot B^s_{p,q} (A_N)}
\leq 
C 
\| f\|_{\dot B^{s+1}_{p,q} (A_D)}. 
\]
The inequalities for $\partial _{x_k} f$ ($k = 1,2,\cdots , n-1$) 
are proved by following the proof of Lemma~\ref{lem_1102-2} 
instead of Lemma~\ref{lem_1102}. 
\end{pf}

\section{Proof of Theorem~\ref{thm_derivative}}

\noindent
{\bf Proof of the well-definedness of $\partial_{x_n}$ in (i) and (iii). } 
Observe that for $M = 0,1,2,\cdots$ 
\[
\| A_N^{ M/2} \partial _{x_n} f \|_{L^1} 
\leq C \| \partial_{x_n} f \|_{\dot B^{ M}_{1,1} (A_N)} 
\leq C \| f \|_{\dot B^{M+1}_{1,1} (A_D)} 
\leq C p_{A_D , M+2} (f),
\]
\[
\| A_N^{\pm M/2} \partial _{x_n} f \|_{L^1} 
\leq C \| f \|_{\dot B^{\pm M+1}_{1,1} (A_D)} 
\leq C q_{A_D , M+2} (f), 
\]
which are assured by the embedding $\dot B^0_{1,1} (A_D) \to L^1 (\Omega)$ 
and Lemma~\ref{lem_1115-1}, it follows that 
\[
p_{A_N,M}(\partial_{x_n}f) \leq C p_{A_D,M+2} (f) \text{ for } f\in \mathcal X(A_D),
\]
\[
q_{A_N,M}(\partial_{x_n}f) \leq C q_{A_D,M+2} (f) \text{ for } f\in \mathcal Z(A_D). 
\]
These give $\partial_{x_n}$ defining maps from 
$\mathcal X(A_D)$, $\mathcal Z(A_D)$ to $\mathcal X(A_N)$, $\mathcal Z(A_N)$, 
respectively. The same argument implies the well-definedness of 
$\partial_{x_n}$ by replacing $A_D$, $A_N$ with each other. 
In the space of distributions, $\partial_{x_n}$ is also well-defined, 
since it is defined by the duality argument. 

\vskip3mm 

\noindent 
{\bf Proof of the boundedness in (ii) and (iii). } 
The result for the Sobolev spaces with $s \geq 1$ 
is obtained by Lemma~\ref{lem_1102}. 
If $s \leq 0$, we regard $\partial_{x_n}$ as a dual operator such that 
\[
_{\dot H^{s-1}_p (A_N)}\langle \partial _{x_n} f, g \rangle _{\dot H^{-s+1}_{p'}(A_N)}
:= 
- _{\dot H^s_p (A_D)} \langle f , \partial _{x_n} g \rangle_{\dot H^{-s}_{p'} (A_D)} . 
\]
We have from Lemma~\ref{lem_1102} that 
\[
| 
_{\dot H^{s-1}_p (A_N)}\langle \partial _{x_n} f, g \rangle _{\dot H^{-s+1}_{p'}(A_N)}
|
\leq \| f \|_{\dot H^s_p (A_D)} \| \partial_{x_n} g\|_{\dot H^{-s}_{p'}A_N)}
\leq C \| f \|_{\dot H^s_p (A_D)} \| g\|_{\dot H^{-s+1}_{p'}A_N)},
\]
which proves that 
\[
\| \partial _{x_n} f \|_{\dot H^{s-1}_p (A_N)} 
\leq C \| f \|_{\dot H ^s_p (A_D)}. 
\]
The case $0 < s < 1$ follows from the complex interpolation of the obtained 
result $s = 0$ and $s = 1$. 
The inhomogeneous case of Sobolev spaces follows similarly. 
The inequality in the Besov spaces are proved by the 
real interpolation of the Sobolev spaces, 
and hence we obtained (ii). 
The boundedness in (iii) follows analogously. 
\hfill $\Box$ 

\vskip3mm 

\noindent 
{\bf Proof of (iv). } 
It is possible to prove (iv) by following the argument for (i), (ii), (iii) 
with Lemma~\ref{lem_1115-1}, Lemma~\ref{lem_1102-2} instead of Lemma~\ref{lem_1102}. 
\hfill $\Box$

\section{Proof of bilinear and trilinear estimates in theorems}

\vskip3mm

\noindent 
{\bf Proof of the Dirichlet case \eqref{1014-1} of Theorem~\ref{thm_bilinear}. }  
Let us start by the case $2 \leq s < 2+1/p$. 
Suppose $f \in \dot H^s_{p_{1}}(A_D) \cap L^{p_3}(\mathbb R^n_+)$, 
$g \in L^{p_2} (\mathbb R^n_+) \cap \dot H^s_{p_2} (A_D)$. 
Lemma~\ref{lem_1} gives that 
\begin{equation}\label{1015-2}
\| A_D ^{s/2} (fg) \|_{L^p} 
\leq C \| \Lambda ^s (fg)_{odd} \|_{L^p (\mathbb R^n)}
= C \left\| \Lambda ^{s-2} (-\Delta) 
      \big( (\sign x_n) f_{odd} \cdot g_{odd}  \big) 
    \right\|_{L^p (\mathbb R^n)}. 
\end{equation}
Here we need to approximate $f_{odd}, g_{odd}$ by smooth odd functions 
to handle their values on $\{x_n = 0\}$.  Put 
\[
F_m := \sum _{ j \leq m} \phi_j (\sqrt{-\Delta}) f_{odd}, 
\quad 
G_m := \sum _{ j \leq m} \phi_j (\sqrt{-\Delta}) g_{odd}, 
\quad m = 1,2,\cdots. 
\]
It is easy to check that $F_m, G_m$ are smooth and 
odd with respect to $x_n$. 
We can see that 
\begin{equation}\label{1015-1}
(-\Delta) \big( (\sign x_n) F_m \cdot G_m  \big) (x)
= 
(\sign x_n)(-\Delta) \big( F_m \cdot G_m  \big) 
\text{ in } \mathcal S'(\mathbb R^n) .
\end{equation}
In fact, for any $\varphi \in \mathcal S'(\mathbb R^n)$ 
\[
\int_{\mathbb R^n} (\sign x_n) F_m G_m \cdot 
(-\Delta) \varphi \, dx 
= I_+ + I_-,  
\quad \text{where }
I_{\pm} := 
\pm \int_{\mathbb R^n_{\pm}}  F_m G_m \cdot 
(-\Delta) \varphi \, dx . 
\]
We have that 
\[
\begin{split}
I_{\pm} 
= 
& 
\pm \Big\{ 
- 
\int_{\partial \mathbb R^n_{\pm}} 
F_m G_m \nabla \varphi \cdot \nu  \, dS
+ 
\int_{\partial \mathbb R^n_{\pm}} 
\nabla \big( F_m G_m \big) \cdot \nu \varphi  \, dS
\\
& 
+ 
\int_{\mathbb R^n_{\pm}} 
\Big( -\Delta \big( F_m G_m \big) \Big) 
\varphi \, dx 
\Big\} . 
\end{split}
\]
The first two terms of $I_\pm$ are zero by 
$F_m G_m, \nabla (F_m G_m) = 0$ on the boundary $\{ x_n = 0 \}$, 
which proves \eqref{1015-1}. 
It follows from \eqref{1015-1} and \eqref{1011-2-0} that 
\[
\begin{split}
\left\| \Lambda ^{s-2} (-\Delta) 
      \big( (\sign x_n) F_m G_m \big) 
    \right\|_{L^p (\mathbb R^n)}
\leq 
& 
\left\| \Lambda ^{s-2} (\sign x_n) (-\Delta) 
      \big( F_m G_m \big) 
    \right\|_{L^p (\mathbb R^n)}
\\
\leq 
&
C  \left\| \Lambda ^{s-2} (-\Delta) \big( F_m G_m \big) 
    \right\|_{L^p (\mathbb R^n)}
\\
= 
& C  \left\| \Lambda ^{s}\big( F_m G_m \big) 
    \right\|_{L^p (\mathbb R^n)}.
\end{split}
\]
The bilinear estimates in the Sobolev spaces in $\mathbb R^n$ gives that 
\[
\begin{split}
\left\| \Lambda ^{s} \big( F_m G_m \big) 
    \right\|_{L^p (\mathbb R^n)}
\leq 
& 
C \big( 
   \| F_m \|_{\dot H^s_{p_1} (\mathbb R^n)} \| G_m \|_{L^{p_2}(\mathbb R^n)} 
   + 
   \| F_m \|_{L^{p_3}(\mathbb R^n)} \| G_m \|_{\dot H^s_{p_4} (\mathbb R^n)}
   \big). 
\end{split}
\]
By taking the limit as $m \to \infty$, we get 
\[
\begin{split}
\left\| \Lambda ^{s} \big( f_{odd} g_{odd} \big) 
    \right\|_{L^p (\mathbb R^n)}
\leq 
& 
C \big( 
   \| f_{odd} \|_{\dot H^s_{p_1} (\mathbb R^n)} \| g_{odd} \|_{L^{p_2}(\mathbb R^n)} 
   + 
   \| f_{odd} \|_{L^{p_3}(\mathbb R^n)} \| g_{odd} \|_{\dot H^s_{p_4} (\mathbb R^n)}
   \big),
\end{split}
\]
where the above convergence is justifiled by the classical theory in 
the whole space case. 
By applying the above inequality and Lemma~\ref{lem_1}, we obtain 
the required estimate \eqref{1014-1}.

We turn to prove the case when $0 < s < 2$ by applying the complex interpolation. 
Lemma~\ref{lem_1} and Bony's paraproduct formula~\cite{Bo-1981} give that 
\begin{equation}\label{1022-3}
\begin{split}
\| A_D ^{s/2} (fg) \|_{L^p} 
\leq \| \Lambda ^s (\sign x_n)(f_{odd} g_{odd}) \|_{L^p (\mathbb R^n)}
\leq II_1 (s) + II_2 (s), 
\end{split}
\end{equation}
where 
\begin{gather}\notag 
II_1 (s):= 
\Big\| \Lambda ^s (\sign x_n) 
   \sum_{ k \geq l+3} 
     \big( \phi_k (\sqrt{-\Delta}) f_{odd} \big) 
     \big( \phi_l (\sqrt{-\Delta}) g_{odd} \big)
\Big\| _{L^p},
\\ \notag 
II_2 (s):= 
\Big\| \Lambda ^s (\sign x_n) 
   \sum_{ k < l+3} 
     \big( \phi_k (\sqrt{-\Delta}) f_{odd} \big) 
     \big( \phi_l (\sqrt{-\Delta}) g_{odd} \big)
\Big\| _{L^p}. 
\end{gather}
Let $\theta$ be such that $s = (1-\theta) \cdot 0 + \theta \cdot 2$. 
The H\"older inequality, the result for the regularity of $s=2$ case and 
the bilinear estimate in the Sobolev spaces in $\mathbb R^n$ imply that 
\[
II_1 (s=0) 
\leq C \| f_{odd} \|_{L^{p_1}(\mathbb R^n)} \| g_{odd} \|_{L^{p_2} (\mathbb R^n)}, 
\quad 
II_1 (s=2) 
\leq C \| f_{odd} \|_{\dot H^2_{p_1} (\mathbb R^n)} 
       \| g_{odd} \|_{L^{p_2} (\mathbb R^n)}.
\]
It follows from the above two inequalities and 
the complex interpolation (see e.g. \cite{BL_1976,Triebel_1978,Triebel_1983})
$(L^p (\mathbb R^n) , \dot H^2 _p (\mathbb R^n)) _\theta
= \dot H^s_p (\mathbb R^n)$ 
that 
\begin{equation}\label{1115-3}
II_1 (s) 
\leq C \| f_{odd} \|_{\dot H^s_{p_1} (\mathbb R^n)} \| g_{odd} \|_{L^{p_2}(\mathbb R^n)}, 
\quad 0 < s < 2.
\end{equation}
Similarly, 
\[
II_2 (s) 
\leq C \| f_{odd} \|_{L^{p_3}(\mathbb R^n)}  \| g_{odd} \|_{\dot H^s_{p_4}(\mathbb R^n)},
\]
which proves \eqref{1014-1} for $0<s<2$, $A = A_D$. 
The Neumann Laplacian case $A = A_N$ for $0 < s < 2$ follows analogously. 
\hfill $\Box$

\vskip3mm 

\noindent 
{\bf Proof of the Neumann case \eqref{1014-2} of Theorem~\ref{thm_bilinear}. }  
We obtain that 
\[
\| A_D ^{s/2} (fg) \|_{L^p} 
\leq C \| \Lambda ^s (f_{even} g_{even}) \|_{L^p (\mathbb R^n)}.
\]
The bilinear estimates in $\mathbb R^n$ give that 
\[
\begin{split}
\| \Lambda ^s (f_{even} g_{even}) \|_{L^p (\mathbb R^n)}
\leq 
&
C ( \| \Lambda ^s f_{even}\|_{L^{p_1} (\mathbb R^n)} 
         \| g_{even} \|_{L^{p_2} (\mathbb R^n)}
       + \| f_{even} \|_{L^{p_3}(\mathbb R^n)} 
         \| \Lambda ^s g_{even} \|_{L^{p_4}(\mathbb R^n)}
   ) 
\\
\leq 
&
C ( \| f\|_{\dot H^s_{p_1} (A_N)} 
         \| g \|_{L^{p_2}}
       + \| f \|_{L^{p_3}} 
         \| g \|_{\dot H^s_{p_4}(A_N)} 
  ), 
\end{split}
\]
which proves \eqref{1014-2}. 
\hfill $\Box$

\vskip3mm 

\noindent
{\bf Remark. } 
There arise no problems for the Neumann case 
such as $\sign x_n$ in contrast to \eqref{1015-2}, 
since $-\Delta (f_{even} g_{even}) = (-\Delta (fg))_{even}$, which is 
observed by that for any sufficiently smooth 
$f_{even}, g_{even}$ and $\varphi \in \mathcal S(\mathbb R^n)$ 
\[
\int_{\mathbb R^n} f_{even} g_{even} 
(-\Delta) \varphi \, dx 
=: I_+ + I_-,  
\quad \text{with }
I_{\pm} := 
\int_{\mathbb R^n_{\pm}} f_{even} g_{even} 
(-\Delta) \varphi \, dx , 
\]
and 
\[
\begin{split}
I_{\pm} 
= 
& 
- 
\int_{\partial \mathbb R^n_{\pm}} 
f_{even} g_{even} \nabla \varphi \cdot \nu  \, dS
+ 
\int_{\partial \mathbb R^n_{\pm}} 
\nabla \big( f_{even} g_{even} \big) \cdot \nu \varphi  \, dS
\\
& 
+ 
\int_{\mathbb R^n_{\pm}} 
\Big( -\Delta \big( \widetilde F_m \widetilde G_m \big) \Big) 
\varphi \, dx .
\end{split}
\]
The sum of the first terms of $I_{\pm}$ is zero by evenness of $f_{even} g_{even} $  
and the second terms of $I_\pm$ are zero by 
oddness of $\partial _{x_n} (f_{even} g_{even} )$ giving 
the well-definedness the value zero on $\{x_n = 0\}$.

\vskip3mm 

\noindent 
{\bf Proof of Theorem~\ref{thm_bilinear_2} for the Besov spaces}. 
Let us start by the Dirichlet Laplacian case. 
We consider a weaker inequality with the Sobolev spaces 
which will be extended to the Besov spaces by means of the real interpolation. 
\[
\| A_D ^{s_0/2}fg \|_{L^p} 
= \| \Lambda ^{s_0 } \sign (x_n) f_{odd} g_{odd} \|_{L^p}
\]

We will apply that the real interpolation of the Sobolev spaces becomes 
the Besov spaces (see~\cite{Iw-2018}) 
and the frequency decomposition such as \eqref{1022-3}. 
Let $0 < s < s_0 < 2 + 1/p$. Then there exists $\theta \in (0,1)$ such that 
\[
\dot B^s_{p,q} (A_D) = (L^p (\mathbb R^n_+) , \dot H^{s_0}_{p} (A_D))_{\theta,q}.
\]
It follows that 
\[
\begin{split}
\| fg \|_{\dot B^s_{p,q}(A_D)} 
\leq  \Big\{ \int_0^\infty \big( t^{-\theta} K\big(t,(f_{odd}g_{odd})_I \big) \big) ^q 
             \, \frac{dt}{t}
      \Big\}^{\frac{1}{q}}
    + 
     \Big\{ \int_0^\infty \big( t^{-\theta} K\big(t,(f_{odd}g_{odd})_{II} \big) \big) ^q 
             \, \frac{dt}{t}
      \Big\}^{\frac{1}{q}},
\end{split}
\]
where 
\[
K(t,a) := 
 \inf 
 \big\{ \| a_0 \|_{L^p(\mathbb R^n)} + t \| a_1 \|_{\dot H^{s_0}_p(\mathbb R^n)} 
       \,\big| \, 
       a = a_0 + a_1, \, a_0 \in L^p (\mathbb R^n) , \, a_1 \in \dot H^{s_0}_{p}(A_D)
 \big\},  
\]
\begin{equation}\label{1115-4}
(f_{odd}g_{odd})_I 
=   \sum_{ k \geq l+3} 
     \big( \phi_k (\sqrt{-\Delta}) f_{odd} \big) 
     \big( \phi_l (\sqrt{-\Delta}) g_{odd} \big),
\end{equation}
\begin{equation}\label{1115-5}
(f_{odd}g_{odd})_{II}
=   \sum_{ k < l+3} 
     \big( \phi_k (\sqrt{-\Delta}) f_{odd} \big) 
     \big( \phi_l (\sqrt{-\Delta}) g_{odd} \big) . 
\end{equation}
Then the bilinear estimates in the Sobolev spaces and the real interpolation 
give that 
\[
\begin{split}
 \Big\{ \int_0^\infty \big( t^{-\theta} K\big(t,(f_{odd}g_{odd})_I \big) \big) ^q 
             \, \frac{dt}{t}
      \Big\}^{\frac{1}{q}}
\leq 
& 
C \| f_{odd} \|_{(L^{p_1}(\mathbb R^n , \dot H^{s_0}_p (\mathbb R^n))_{\theta ,q}} 
       \| g_{odd} \|_{L^{p_2} (\mathbb R^n)}
\\
\leq 
&  
C \| f \|_{\dot B^s_{p_1,q}} \| g \|_{L^{p_2}}, 
\end{split}
\]
\[
\begin{split}
 \Big\{ \int_0^\infty \big( t^{-\theta} K\big(t,(f_{odd}g_{odd})_{II} \big) \big) ^q 
             \, \frac{dt}{t}
      \Big\}^{\frac{1}{q}}
\leq 
& 
C \| f_{odd} \|_{L^{p_3} (\mathbb R^n)}
  \| g_{odd} \|_{(L^{p_4}(\mathbb R^n , \dot H^{s_0}_{p_4} (\mathbb R^n))_{\theta ,q}} 
\\
\leq 
&  
C \| f \|_{L^{p_3}} \| g \|_{\dot B^s_{p_4 , q}}, 
\end{split}
\]
which proves the result for the homogeneous Besov spaces for $1 < p < \infty$. 
The bilinear estimates for the inhomogeneous Besov spaces also follows 
from the those in homogeneous ones.

The case $p = 1, \infty$ needs some modification. 
If $p = 1$, we take $s$ satisfying $2 < s < 2 + 1/p =3$. 
Observe that for $0 < \alpha < 1/p = 1$ 
\[
\| \Lambda ^{\alpha} \big( (\sign x_n) F \big) \|_{L^1} \leq 
C \| F \|_{\dot B^{\alpha}_{1,1} (\mathbb R^n)}, 
\quad F \in \dot B^{\alpha}_{1,1} (\mathbb R^n),
\]
which is proved analogously to \eqref{1011-2-0} 
and by applying $\dot B^0_{1,1} (\mathbb R^n) \to L^1 (\mathbb R^n)$ 
and the boundedness of the Fourier multiplier $|\xi_n|^{\alpha}/|\xi|^{\alpha}$ 
in $\dot B^0_{1,1} (\mathbb R^n)$ to \eqref{1115-1}. Then 
the inequality with $L^1$ norm in the left hand side replaced by the Besov norm 
$\dot B^0_{1,q}(A_D)$ also hold thanks to the real interpolation, 
and we have 
\[
\begin{split}
\| fg \|_{\dot B^s_{1,q} (A_D)} 
\leq 
& 
C \| (\sign x_n) (-\Delta) f_{odd} g_{odd}  \|_{\dot B^{s-2}_{1,q}(\mathbb R^n)}
\leq C \| (-\Delta) f_{odd} g_{odd}  \|_{\dot B^{s-2}_{1,q} (\mathbb R^n)}
\\
=
& 
 C \| f_{odd} g_{odd} \|_{\dot B^s_{1,q}(\mathbb R^n)}. 
\end{split}
\]
The classical bilinear estimates in $\mathbb R^n$ and Lemma~\ref{lem_1115-1} 
give that 
\[
\| f_{odd} g_{odd} \|_{\dot B^s_{1,q}(\mathbb R^n)} 
\leq C \big( \| f \|_{\dot B^s_{p_1} (A_D)} \| g \|_{L^{p_2}} 
    + \| f \|_{L^{p_3}} \| g \|_{\dot B^s_{p_4 ,q} (A_D)} 
    \big),
\]
which proves the bilinear estimate of homogeneous spaces for $p = 1$ 
and $2 < s < 3$. The case when $0 < s \leq 2$  follows from the frequency 
decomposition $f_{odd} g_{odd} = (f_{odd} g_{odd})_I + (f_{odd} g_{odd})_{II}$ 
and the real interpolation as previous \eqref{1115-3}. 
As for the case when $ p = \infty$, we start by 
\[
\| A_D (fg) \|_{L^\infty} 
\leq C \|\Delta (f_{odd} g_{odd}) \|_{L^\infty (\mathbb R^n)}, 
\]
where we applied Lemma~\ref{lem_1}, \eqref{1015-1}. By decomposing as 
\eqref{1115-4}, \eqref{1115-5} and the bilinear estimates in $\mathbb R^n$, 
we have 
\[
\|\Delta (f_{odd} g_{odd})_I \|_{L^\infty (\mathbb R^n)} 
\leq C \| f_{odd} \|_{\dot B^2_{\infty,1}(\mathbb R^n)} 
       \| g_{odd} \|_{L^\infty (\mathbb R^n)},
\]
\[
\|\Delta (f_{odd} g_{odd})_{II} \|_{L^\infty (\mathbb R^n)} 
\leq C \| f_{odd} \|_{L^\infty(\mathbb R^n)} 
       \| g_{odd} \|_{\dot B^2_{\infty,1} (\mathbb R^n)}. 
\]
These and the real interpolation imply that for any $0 < s < 2$
\[
\| (f_{odd} g_{odd})_I \|_{\dot B^s_{\infty,q}(\mathbb R^n)} 
\leq C \| f_{odd} \|_{\dot B^s_{\infty,q}(\mathbb R^n)} 
       \| g_{odd} \|_{L^\infty (\mathbb R^n)},
\]
\[
\| (f_{odd} g_{odd})_{II} \|_{\dot B^s_{\infty,q}(\mathbb R^n)} 
\leq C \| f_{odd} \|_{L^\infty(\mathbb R^n)} 
       \| g_{odd} \|_{\dot B^s_{\infty,q} (\mathbb R^n)},
\]
which prove the bilinear estimate in $\dot B^s_{\infty,q} (A_D)$. 
The estimates for the inhomogeneous Besov spaces follows analogously. 

The Neumann Laplacian case is proved by following the above argument, 
and we notice that $s_0$ can be choosen as arbitrary positive number 
as well as the Sobolev spaces. 
\hfill $\Box$

\vskip3mm 

\noindent 
{\bf Proof of the trilinear estimates in Corollary~\ref{cor_tri}. }
Observing that odd extention of $fgh$ is given by $f_{odd} g_{odd} h_{odd}$, 
we obtain that 
\begin{equation}\label{1101-1}
\| A_D ^{s/2} (fgh) \|_{L^p} 
\leq C \| \Lambda ^s (f_{odd} g_{odd} h_{odd}) \|_{L^p (\mathbb R^n)}.
\end{equation}
The trilinear estimate in $\mathbb R^n$ gives the results. 
\hfill $\Box$

\section{Counter examples in Theorem \ref{thm_counter_example} 
and Theorem~\ref{thm_bilinear_2}}

\noindent 
\underline{The case when $n = 1$}.  
We construct $f,g $ such that 
$f,g \in H^s_p(A_D)$ for any $s \geq 0$ and $1\leq p\leq \infty$ but 
$fg \not \in H^{2+\frac{1}{p}}_{p} (A_D)$. 
Let $\varphi$ be such that 
\begin{equation}\label{1022-2}
\varphi \in C^\infty ([0,\infty)), \quad 0 \leq \varphi \leq 1, 
\quad 
\varphi (x)= \begin{cases}
1 & \text{for } 0 \leq x \leq 1/2, 
\\
0 & \text{for } x \geq 1.
\end{cases}
\end{equation}
Take $f,g$ such that 
\[
f(x) = g(x) = x \varphi (x).
\]
It is easy to show that $f, g \in H^s_p (A_D)$ for any $s \geq 0$, 
$1\leq p \leq \infty$. 
It suffices to prove $A_D (fg) \not \in H^{1/p}_p (A_D)$. 
We see that 
\[
A_D (fg) = 
(A_D f ) g - \partial _x f \cdot \partial_x g + f A_D g.
\]
The first  and the third term are in 
$\dot H^\frac{1}{p}_p (A_D)$, since they are in $C_0 ^\infty (0,\infty)$. 
The second term is 
\[
\partial _x f \cdot \partial_x g
=
\varphi^2 + 2 x \varphi \varphi ' + x^2 (\varphi') ^2 . 
\]
Since $2 x \varphi \varphi ',  x^2 (\varphi') ^2 \in C_0 ^\infty (0,\infty)$, 
they belong to $H^{1/p}_p (A_D)$. Put 
\[
\Phi := \varphi ^2.
\]
Noting that $A_D ^{1/2p} \Phi \in L^p (\mathbb R_+)$ is equivalent to 
$\Lambda ^{1/p} \Phi_{odd} \in L^p (\mathbb R)$, we should consider 
\[
\Lambda ^{1/p} \Phi_{odd} (x) 
= C \int_{\mathbb R} 
  \frac{\Phi_{odd} (x) - \Phi_{odd} (y)}{|x-y|^{1+\frac{1}{p}}}
  \, dy ,
\]
and one can see that there exist $c > 0$ and 
$\delta > 0$ such that 
\[
\Lambda ^{1/p} \Phi_{odd} (x) 
\geq \frac{c}{|x|^{\frac{1}{p}}} \text{ if } 0 < x < \delta, 
\quad 
\Lambda ^{1/p} \Phi_{odd} (x) 
\leq \frac{-c}{|x|^{\frac{1}{p}}} \text{ if } -\delta < x < 0 .
\]
Hence, we get $\Lambda^{1/p} \Phi_{odd} \not \in L^p (\mathbb R)$, 
which proves that $fg \not \in \dot H^{2+\frac{1}{p}}_p (A_D)$.

As for counter example in the Besov spaces 
$\dot B^{2+\frac{1}{p}}_{p,q} (A_D)$ with $1 \leq q < \infty$, 
we can also prove that 
$\Phi _{odd} \not \in \dot B^{\frac{1}{p}}_{p,q} (\mathbb R) $ 
for $1 \leq q < \infty $. In fact, it follows that
\[
\begin{split}
2^{\frac{1}{p}j} 
\| \phi_j (\sqrt{-\Delta}) \Phi_{odd} \|_{L^p(\mathbb R)}
\geq 
& 
\| \phi_j (\sqrt{-\Delta}) \Phi_{odd} \|_{L^\infty(\mathbb R)} .
\end{split}
\]
We have 
\[
\begin{split}
\phi_j (\sqrt{-\Delta}) \Phi_{odd} (2^{-j} x )
= 
& 2^j \int _{\mathbb R} 
  \phi_0 (|x - 2^j y|) \Phi_{odd} (y) dy 
\\
= 
& \int _{\mathbb R} 
  \phi_0 (|x - y|) \Phi_{odd} ( 2^{-j}y) dy 
\\
= 
& \int _{0}^\infty 
  \big( \phi_0 (|x - y|) - \phi_0(|x+y|) \Big) \varphi  ( 2^{-j}y) ^2 dy 
\\
\to 
& 
\int _{0}^\infty 
  \big( \phi_0 (|x - y|) - \phi_0(|x+y|) \Big) dy 
 ( \not = 0) 
 \text{ as } j \to \infty ,
\end{split}
\]
which proves that for some $j_0 \in \mathbb Z$ and $c > 0$ 
\[
\Big\{ 
\sum _{j \in \mathbb Z} 
  \Big( 2^{\frac{1}{p}j} 
        \| \phi_j (\sqrt{-\Delta}) \Phi_{odd} \|_{L^p(\mathbb R)}
   \Big) ^q
\Big\} ^{\frac{1}{q}} 
\geq 
\Big\{ 
\sum _{j \geq j_0} c^q
\Big\} ^{\frac{1}{q}} 
= \infty \quad \text{ if } q < \infty . 
\]

\vskip3mm 

\noindent 
\underline{The case when $n \geq 2$}. For $\varphi$ satisfying \eqref{1022-2}, 
put 
\[
f = g= x_n \varphi (x_n) \cdot \Big( \varphi (x_1) \cdots \varphi (x_{n-1})\Big) . 
\]
We should consider 
\[
A_D (f g) = (A_D f) g - \nabla f \cdot \nabla g + f A_D g, 
\]
and the terms except for $\partial_{x_n} f \partial _{x_n g}$ 
are in $\dot H^{1/p} _p (\mathbb R^n)$ but the second term is 
\[
-\partial_{x_n} f \partial _{x_n g}
= 
-\Big( 
\varphi (x_n)^2 + 2 x_n \varphi(x_n) \varphi ' (x_n) + x_n^2 (\varphi'(x_n))^2
\Big) 
( \varphi (x_1)^2 \cdots \varphi (x_{n-1}) ^2).
\]
Similarly to the case when $n = 1$, 
the above terms having derivative $\varphi '$ in the right hand side 
is in $C^\infty _0 (\mathbb R^n_+)$, but 
for the first one 
$\widetilde \Phi := \varphi (x_n)^2 \varphi (x_1)^2 \cdots \varphi (x_{n-1}) ^2$, 
we can show that there exist $c > 0$ and $\delta > 0$ such that 
for $|x| \leq \delta $ 
\[
\Lambda ^{1/p} \widetilde \Phi_{odd} (x) 
\geq \frac{c}{|x_n|^{\frac{1}{p}}} \text{ if } 0 < x_n < \delta, 
\quad 
\Lambda ^{1/p} \widetilde \Phi_{odd} (x) 
\leq \frac{-c}{|x_n|^{\frac{1}{p}}} \text{ if } -\delta < x < 0 
\]
which proves that $\Lambda ^{1/p}\widetilde\Phi _{odd} \not \in L^p (\mathbb R^n) $. 
Therefore $fg \not \in \dot H^{2+\frac{1}{p}}_p (A_D)$.

\section{Sobolev spaces}\label{sec_Sobolev}

In this section, let us explain that we can verify the well-definedness of 
the Sobolev spaces $H^s_p (A)$ and $\dot H^s_p (A)$ for $A = A_D, A_N$.

\begin{prop}
Let $A = A_D$ or $A_N$, $s \in \mathbb R$, $1 \leq p \leq \infty$. 
\begin{enumerate}
\item[(i)] $H^s_p (A)$, $\dot H^s_p (A)$ are Banach spaces, and enjoy 
\[
\mathcal X(A) \hookrightarrow H^s_p(A) \hookrightarrow \mathcal X'(A), 
\quad 
\mathcal Z(A) \hookrightarrow \dot H^s_p(A) \hookrightarrow \mathcal Z'(A). 
\]

\item[(ii)] Let $1 \leq p < \infty$ and $1/p + 1/p' = 1$. 
Then the dual spaces of $H^s_p (A)$, $\dot H^s_p(A)$ are 
$H^{-s}_{p'} (A)$, $\dot H^{-s}_{p'} (A)$, respectively. 

\item[(iii)] Let $\alpha \in \mathbb R$. Then 
\[
(1+A)^{\alpha /2  } f \in H^{s}_p (A) 
\text{ for } f \in H^{s+\alpha}_p(A), 
\quad 
A^{\alpha /2  } f \in \dot H^{s}_p (A) 
\text{ for } f \in \dot H^{s+\alpha}_p(A)
\]

\item[(iv)] 
Let $1 < p \leq r < \infty$. Then 
\[
H^{s+n(\frac{1}{p}-\frac{1}{r})}_p (A) 
\hookrightarrow H^s_r(A), 
\quad 
\dot H^{s+n(\frac{1}{p}-\frac{1}{r})}_p (A) 
\hookrightarrow \dot H^s_r(A). 
\]

\item[(v)] Let $s < n/p$. Then 
\[
\dot H^s_p (A) \simeq 
\Big\{ f \in \mathcal X'(A) \, \Big| \, 
 f = \sum _{j \in \mathbb Z} \phi_j (\sqrt{A})f 
  \text{ in } \mathcal X'(A) , 
  \,\, \| f \|_{\dot H^s_p(A)}  < \infty 
\Big\} . 
\]
\end{enumerate}

\end{prop}

\begin{pf} 
Let us prove for the homogeneous spaces only, since the inhomogeneous 
case follows analogously with a modification of the proof below 
by replacing $\mathcal Z$, $\mathcal Z'$, the operator $A^{s/2}$ 
with $\mathcal X$, $\mathcal X'$ , the operator $(1+A)^{s/2}$,  respectively. 

Step 1. 
It is sufficient to show the completeness to prove the spaces are Banach spaces. 
Let $\{ f_N \}_{N=1}^\infty $ be a Cauchy sequence in $\dot H^s_p (A)$. 
Then $\{ A^{s/2} f_N \}_{N=1}^\infty$ is a Cauchy sequence in 
$L^p (\mathbb R^n_+)$, whose completeness gives that 
$F \in L^p (\mathbb R^n_+)$ exists such that $A^{s/2} f_N$ converges to 
$F$ in $L^p (\mathbb R^n_+)$ as $N \to \infty$. 
Let $f$ be a element of $\mathcal Z'(A)$ given by 
\[
f:= A^{-s/2} F \text{ in } \mathcal Z'(A),
\]
where we note that the well-definedness of 
$A^{s/2} : \mathcal Z'(A) \to \mathcal Z'(A)$ 
is already known in the paper~\cite{IMT-preprint2} (see also \cite{Iw-2018}). 
Then we find that $f_N$ tends to $f$ in $\dot H^s_p(A)$ as $N \to \infty$. 
As for the continuous embedding, 
for $M\in \mathbb N$ with $M > s + n(1-1/p)$ and $f \in Z(A)$, 
\[
\begin{split}
\| f \|_{\dot H^s_p (A)} 
\leq 
& C \sum _{j \in \mathbb Z} 2^{sj} \| \phi_j (\sqrt{A}) f\|_{L^p} 
\leq 
 C \sum _{j \in \mathbb Z} 2^{sj+n(1-\frac{1}{p})j} \| \phi_j (\sqrt{A}) f\|_{L^1} 
\\
\leq 
& C \Big( \sum _{j \in \mathbb Z} 2^{sj+n(1-\frac{1}{p})j -M|j|} \Big) q_{A,M} (f),
\end{split}
\]
which proves $\mathcal Z(A) \hookrightarrow \dot H^s_p (A)$. 
The second embedding is verified by 
\[
|_{\mathcal Z'}\langle f,g \rangle_{\mathcal Z} |
= |_{\mathcal Z'}\langle A^{s/2}f, A^{-s/2}g \rangle_{\mathcal Z} |
\leq \| A^{s/2}f \|_{L^p} \| A^{-s/2}g \|_{L^{p'}} 
\leq C \| f \|_{\dot H^s_p (A)} q_{A,M'} (g),
\]
where $g \in \mathcal Z(A)$, $M' \in \mathbb N$ satisfies 
$M' > -s + n (1-1/p')$ with $1/p + 1/p' = 1$. 

Step 2. Let us prove the duality. For $f \in \dot H^{-s}_{p'}(A)$, let 
$T_f$ be defined by 
\[
T_f (g) := \int_{\mathbb R^n_+} (A^{-s/2} f)  \overline{A^{s/2} g} \, dx, 
\quad g \in \dot H^s_p (A).
\]
Then we have $\dot H^{-s}_{p'} (A) \hookrightarrow (\dot H^s_p(A))'$ by 
\[
|T_f (g)| \leq \| f \|_{\dot H^{-s}_{p'}(A)} \| g\|_{\dot H^s_p (A)} .
\]

Conversely, Let $F \in (\dot H^s_p(A))'$ and define 
\[
T(G) := F( A^{-s/2} G), \quad G \in L^p (\mathbb R^n_+). 
\]
It follows that 
\[
|T(G)| \leq \| F \|_{(\dot H^s_p(A))'} \| A^{-s/2} G \|_{\dot H^s_p (A)}
= \| F \|_{(\dot H^s_p(A))'} \|  G \|_{L^p}.
\]
Since $(L^p (\mathbb R^n_+))' = L^{p'} (\mathbb R^n_+)$, 
$\widetilde f \in L^{p'} (\mathbb R^n_+) $ exists such that 
\[
T(G) = \int_{\mathbb R^n_+} \widetilde f(x) \overline{G(x)} \, dx 
\quad \text{ for any } G \in L^p (\mathbb R^n_+). 
\]
Observe that for any $g \in \dot H^s_p (A)$, 
\[
F(g) = T(A^{s/2} g) 
= \int_{\mathbb R^n_+} \widetilde f(x) \overline{A^{s/2} g (x)} \, dx, 
\]
define $f$ and $\langle f ,g \rangle$ by 
\[
f := A^{s/2} \widetilde f \in \dot H^{-s}_p (A), 
\quad 
\langle f ,g \rangle := \int_{\mathbb R^n_+} f (x) \overline{g(x)} \, dx 
\, \text{ for }g \in \dot H^s_p (A). 
\]
We obtain for any $g \in \dot H^s_p (A)$ 
\[
F(g) = \langle f ,g \rangle, 
\quad 
\| f\| _{\dot H^{-s}_p (A)} 
\leq \| F \|_{(\dot H^s_p (A))'},
\]
which proves $ (\dot H^s_p(A))'\hookrightarrow \dot H^{-s}_{p'} (A)$. 

Step 3. We prove the lifting property in this step. 
Let $f \in \dot H^{s+\alpha}_p (A)$. 
Since $A^{\alpha/2}$ is a operator from $\mathcal Z'(A)$ to itself, 
$A^{\alpha /2} f \in \mathcal Z'(A)$, and the definition of 
$\dot H^s_p(A)$ implies 
\[
\| A^{s/2} f \|_{\dot H^s_p (A) }  \leq 
\| f \|_{\dot H^{s+\alpha}_p} .
\]

Step 4. We prove the embedding theorem in this step. 
Let $\widetilde f$ be $f_{odd}$ for $A = A_D$, $f_{even}$ for $A =A_N$ to 
apply Lemma~\ref{lem_1}. 
If $s \geq 0$, the Sobolev embedding in $\mathbb R^n$ gives that 
\[
\| f \|_{\dot H^s_r (A)} 
\leq \| \Lambda ^s \widetilde f \|_{L^r (\mathbb R^n)}
\leq C \| \Lambda ^{s+ n(\frac{1}{p}-\frac{1}{r})} \widetilde f \|_{L^p(\mathbb R^n)}
\leq C \| f \|_{\dot H^{s+n(\frac{1}{p}-\frac{1}{r})}_p (A)} .
\]
The lifting property obtained in Step 3 proves the case $s < 0$. 

Step 5. We prove the characterization of $\dot H^s_p (A)$ as a subspace of 
$\mathcal X'(A)$ in this step following the argument in some literature 
following the argumetn as in e.g. \cite{IMT-preprint2,KoYa-1994}. 
Let $f \in \dot H^s_p (A)$ where $s < n/p$. The resolution of identity 
in $\mathcal Z'(A)$ (see~\cite{IMT-preprint2}) gives that 
\[
f = \Big( \sum _{ j \leq 0} + \sum _{ j > 0} \Big) \phi_j(\sqrt{A}) f 
\quad \text{in } \mathcal Z'(A).
\]
It is sufficient to justify this expansion in $\mathcal X'(A)$. 
We can see the high spectral component is regarded as an element of 
$\mathcal X'(A)$. For the low spectral component, 
it is sufficient to show that it belongs to 
$L^\infty (\mathbb R^n_+)$, 
which is assured by 
\[
\begin{split}
\Big\| \sum _{ j < 0} \phi_j (\sqrt{A} ) f \Big\|_{L^\infty} 
\leq 
& 
C \sum _{ j < 0} \| \phi_j (\sqrt{A}) f \|_{L^\infty} 
\leq 
C \sum _{ j < 0} 2^{\frac{n}{p}j } \| A^{-s/2} \phi_j (\sqrt{A}) A^{s/2}f \|_{L^p} 
\\
\leq 
& 
C \Big( \sum _{ j < 0} 2^{(-s+\frac{n}{p})j } \Big) \| A^{s/2} f \|_{L^p} . 
\end{split}
\]
Hence we obtained (v) by the embedding 
$L^\infty (\mathbb R^n_+) \hookrightarrow \mathcal Z'(A)$. 
\end{pf}

\vskip3mm 

\noindent
{\bf Acknowledgements.}
The author was supported by the Grant-in-Aid for Young Scientists (A) (No.~17H04824)
from JSPS.

\vskip3mm 
%
%
%
%

\end{document}